\documentclass{article}[15pt]
\usepackage{amsmath,amsmath, epsfig,amssymb, multicol}
\usepackage[active]{srcltx}
\usepackage{epstopdf}
\usepackage{graphicx}
\usepackage{color}
\usepackage{setspace}
\newtheorem{lemma}{Lemma}
\newtheorem{theorem}{Theorem}

\newtheorem{fact}{Fact}

\newcommand{\h}{{\cal  H}}
\newcommand{\f}{{\cal  F}}
\newcommand{\ex}{{\rm ex}}
\newcommand{\EX}{{\rm EX}}
\newcommand{\C}{{\cal C}}

\begin{document}
\title{Tur\'an numbers of cycles plus a general graph}

\author{{Chunyang Dou\footnote{School of Mathematical Sciences, Anhui University,
Hefei  230601, P.~R.~China. Email:{\tt chunyang@stu.ahu.edu.cn.} Supported in part by the National Natural Science Foundation of China
(No.\ 12071002, No.\ 12471319) and the Anhui Provincial Natural Science Foundation (No. 2208085J22).}}
~~~
{Fu-Tao Hu\footnote{ Center for Pure Mathematics, School of Mathematical Sciences, Anhui University,
Hefei  230601, P.~R.~China. Email: {\tt hufu@ahu.edu.cn}. Supported in part by the National Natural Science Foundation of China
(No.\ 12331012).}}
~~~
 {Xing Peng\footnote{Center for Pure Mathematics, School of Mathematical Sciences, Anhui University,
Hefei  230601, P.~R.~China. Email: {\tt x2peng@ahu.edu.cn}. Supported in part by the National Natural Science Foundation of China
(No.\ 12071002, No.\ 12471319) and the Anhui Provincial Natural Science Foundation (No. 2208085J22). } }}
\maketitle

\begin{abstract}
For a family of graphs $\f$, a graph $G$ is $\f$-free if it does not contain a member of $\f$ as a subgraph.
The Tur\'an number $\ex(n,\f)$ is the maximum number of edges in an $n$-vertex graph which is $\f$-free.  Let $\C_{\geq k}$ be the set of cycles with length at least $k$.
In this paper, we investigate  the Tur\'an number of $\{\C_{\geq k}, F\}$ for a general 2-connected graph $F$. To be precise, we determine $\ex(n, \{\C_{\geq k}, F\})$ apart from a constant additive term, where $F$   is either a 2-connected nonbipartite graph or  a 2-connected bipartite graph under some conditions. This is an extension of a previous result on the Tur\'an number of $\{\C_{\geq k}, K_r\}$ by the first author, Ning, and the third author.

%
%

\end{abstract}

\section{Introduction}
This paper studies the Tur\'an number of graphs. Let $\f$ be a family of graphs. A  graph $G$ is said {\it $\f$-free} if it does not contain a member of $\f$ as a subgraph. The {\it Tur\'an number} $\ex(n,\f)$ is the maximum number of edges in an $n$-vertex $\f$-free graph. An $n$-vertex $\f$-free graph with  $\ex(n,\f)$ edges is called an {\it extremal graph}. The set of extremal graphs with $n$ vertices is denoted by $\EX(n,\f)$.
If $\f$ contains only one graph $F$, then we will write $\ex(n,F)$ and $\EX(n,F)$ instead of $\ex(n,\f)$ and $\EX(n,\f)$ respectively.

Perhaps Mantel's result on  $\ex(n,K_3)$ is the earliest one studying the Tur\'an number of a graph in  history. The famous Tur\'an's theorem \cite{Turan} gives us the value of $\ex(n,K_r)$. Tur\'an's theorem is believed as the start of extremal graph theory. In an influential paper, Erd\H{o}s and Gallai \cite{EG59} investigated the Tur\'an number of a matching $M_{s+1}$, a path $P_k$, and the set of cycles with length at least $k$ (denoted by $\C_{\geq k}$).
We only state the result on $\ex(n,\C_{\geq k})$ as it will be needed for us.
\begin{theorem}[Erd\H{o}s-Gallai \cite{EG59}]\label{EGcycle}
$\ex(n, \C_{\geq k}) \leq \frac{(k-1)(n-1)}{2}$, where $k \geq 3$.
\end{theorem}

In the study of Tur\'an number of graphs, the celebrated Erd\H{o}s-Stone-Simonovits Theorem \cite{ES,ES1}  gives an asymptotic formula for $\ex(n,\f)$:
\[
\ex(n, \f) = \left(1 - \frac{1}{c(\f)} \right) \binom{n}{2} + o(n^2),
\]
where $c(\f)=\min\{\chi(F)-1: F \in \f \}$.

 There is a rich literature on the Tur\'an number of single graphs. If one wishes to study the Tur\'an number of $\f$ with two graphs, then it is necessary to assume that $\f$ contains a bipartite graph because of  Erd\H{o}s-Stone-Simonovits Theorem. Chv\'{a}tal and Hanson's  \cite{CH}  classical result in this direction gives us the Tur\'an number of $\{M_{\nu},K_{1,\Delta}\}$. Note that Abbott, Hanson, and Sauer \cite{AHS} previously studied
   a special case where $\nu=\Delta=k$.

   Inspired by a recent work of Alon and Frankl \cite{AF} on the Tur\'an number of $\{M_{s},K_{r}\}$, the study of Tur\'an number of $\f$ with two graphs and its variants  turns to be fruitful. Wang, Hou, and Ma \cite{WHM} established a spectral analogue of the result by Alon and Frankl.
  Following the study in the paper by Erd\H{o}s and Gallai,  Katona and Xiao \cite{KX} determined the value of $\ex(n,\{P_k,K_r\})$ for $r \leq \lfloor k/2 \rfloor$ and sufficiently large $n$. They also conjectured the value of $\ex(n,\{P_k,K_r\})$ for $r > \lfloor k/2 \rfloor$ and  $n$ being large. Fang, Zhu, and Chen \cite{FZC} recently proved this conjecture. Note that the first author, Ning, and the third author \cite{DNP} confirmed this conjecture in a stronger form independently. Additionally, Katona and Xiao \cite{KX} proposed an intriguing conjecture on the Tur\'an number  of $\{P_k, F\}$ for a connected nonbipartite graph $F$. A recent paper by Liu and Kang \cite{LK} verified this conjecture. It is natural to ask for the determination of the Tur\'an number of $\ex(n,\{M_s, F\})$ for a nonbipartite graph $F$. This problem was solved by Zhu and Chen in the paper \cite{ZC}. Actually, they also managed to show the value of $\ex(n,\{M_s, F\})$ for a bipartite $F$ with some assumptions. There are several other related results in this direction, for example \cite{DG,LLK,MHY,XK}.

 It is a well-known fact that $\ex(n,P_k)$ is a corollary of   $\ex(n,\C_{\geq k+1})$. Motivated by the fact above,
 the first author, Ning, and the third author \cite{DNP} investigated the Tur\'an number of  $\ex(n,\{\C_{\geq k+1},K_r\})$. To introduce results from \cite{DNP}, we need the following definition of the join of two graphs.
Let $F_1$ and $F_2$ be two vertex disjoint graphs.  The {\it join}  $F_1 \vee F_2$  is a graph obtained from $F_1 \cup F_2$ by connecting each vertex of $F_1$ and each vertex of $F_2$.

  Assume that $n-1=p(k-2)+q$ with $q \leq k-3$. Let $F(n,k,r)$  be the graph which consists of $p$ copies of $T(k-1,r-1)$ and one copy of $T(q+1,r-1)$ sharing a common vertex.  The number of edges in $F(n,k,r)$ is denoted by $f(n,k,r)$. Note that
\[
f(n,k,r)=p\ex(k-1,K_r)+\ex(q+1,K_r).
\]
 Additionally, let
$$G_1=K_{\lfloor (k-1)/2 \rfloor} \vee I_{n-\lfloor (k-1)/2 \rfloor }$$
and
$$G_2=T(\lfloor (k-1)/2 \rfloor,r-2) \vee I_{n-\lfloor (k-1)/2 \rfloor }.$$
One of the main results from \cite{DNP} is the following one.
\begin{theorem} [Dou, Ning, Peng \cite{DNP}]\label{thm1}
Let  $\f=\{\C_{\ge k},K_r\}$ and $n-1=p(k-2)+q$. If $k\geq 6$ is even and  $\lfloor (k-1)/2 \rfloor+2\le r<k \leq n$, then
\[
\ex(n,\f)=f(n,k,r).
\]
 If $k \geq 5$ is odd and  $\lfloor (k-1)/2 \rfloor+2\le r<k \leq n$, then
\[
\ex(n,\f)=\max\left\{f(n,k,r), e(G_1)\right\}.
\]
\end{theorem}
The authors \cite{DNP} also determined the exact value of $\ex(n,\{P_k,K_r\})$ for all $n$ and $r > \lfloor k/2 \rfloor$, which confirms Katona and Xiao's conjecture in a stronger form.
The following result addresses the case where $r\le\lfloor (k-1)/ 2\rfloor+1$.
\begin{theorem} [Dou, Ning, Peng \cite{DNP}] \label{thm2}
Let  $\f=\{\C_{\ge k},K_r\}$. If $r\le\lfloor (k-1)/ 2\rfloor+1$ and $n \ge k\ge 5$, then
\[
\ex(n,\f)=e(G_2)
\]
provided   $n \geq \tfrac{k^3}{4}$.
\end{theorem}
 As the value of  $\ex(n,\{P_k,F\})$  is known for a connected nonbipartite graph $F$, which motivates us to study the Tur\'an number of  $\{\C_{\geq k},F\}$ in this paper. Before presenting our results, we need to introduce a number of concepts.
 A {\it vertex covering} $S$ of  graph $F$ is a subset of $V(F)$ such that $V(F-S)$ is an independent set.  Let
 \[
 \h=\{F[S]: S \text{ is a vertex covering of } F\}.
 \]
 By the definition, we know $\h \neq \emptyset$. For a bipartite graph $F$, let $p(F)$ be  the smallest possible coloring class among all 2-colorings of $H$. For a nonbipartite graph $F$, let $p(F)=\infty$.
Extending techniques from \cite{DNP}, we manage to determine  $\ex(n,\{\C_{\geq k},F\})$ apart from a constant additive term for large $n$, where $F$ is a 2-connected graph with $p(F) \geq \lfloor(k-1)/2\rfloor+1$.
	Throughout the paper, for  a positive integer $k$, we fix
$$t=\lfloor(k-1)/2\rfloor.$$
To determine $\ex(n,\{\C_{\geq k},F\})$, we follow the approach in \cite{DNP} and consider the following variant.  Let $\ex_{2\textrm{-conn}}(n,\f)$ be the maximum number of edges among all $n$-vertex graphs which are 2-connected  and $\f$-free. We will prove the following result for odd $k$.
\begin{theorem}\label{odd}
	 Let  $\f=\{\C_{\ge k}, F\}$. If $F$ is a  graph with $p(F)\ge t+1$ and $k \geq 5$ is odd, then
	 \[\ex_{2\textrm{-conn}}(n,\f)=\ex(t,\h)+t(n-t)\]
	 for sufficiently large $n$.
	\end{theorem}
Note that if $k=3$, then an $\f$-free graph must be a forest which is not 2-connected.
For even $k$, we need to introduce an additional definition. For a graph $F$, let
\[
\h'=\{F-\{u,v\}: u \textrm{ is adjacent to } v \textrm{ in } F\}.
\]
Our result for even $k$ is the following one.
	\begin{theorem}\label{even}
		Let  $\f=\{\C_{\ge k}, F\}$. If $F$ is a  graph with $p(F)\ge t+1$ and $k \geq 6$ is even, then
		\[
		\ex_{2\textrm{-conn}}(n,\f)=\ex(t,\h)+t(n-t)+\ell
		\]
		for sufficiently large $n$, where $\ell \in \{0,1\}$. If  every graph in $\EX(t,\h)$ contains  a member of $\h'$ as a subgraph, then $\ell=0$.
	\end{theorem}
For $k=4$, if an $\f$-free graph $G$ contains cycles, then cycles are triangles and any two triangles can not share an edge. Thus $G$ can not be 2-connected and  we assume $k \geq 6$.

Let $\f=\{\C_{\geq k}, F\}$ with $p(F) \geq t+1$.
If $k$ is odd, then the combination of the lower bound in Theorem \ref{odd} and Theorem \ref{EGcycle} yields
\[
\ex(n,\f)=\frac{k-1}{2}n+O_k(1)
\]
for sufficiently large $n$. This shows that the case where $k$ is odd is not interesting.  For even $k$, our result is as follows.
\begin{theorem} \label{nocon2}
	Let  $\f=\{\C_{\ge k}, F\}$. If $k$ is even and  $F$ is a 2-connected graph with $p(F)\ge t+1$,  then
	\[
	\ex(n,\f)=n\max\left\{\frac{k-2}{2},\frac{\ex(k-1,F)}{k-2}\right\}+O_k(1)
	\]
for sufficiently large $n$.
\end{theorem}

The organization of this paper is as follows. We will prove a number of necessary lemmas in Section 2. Proofs for Theorem \ref{odd} and Theorem \ref{even} will be presented in Section 3. Section 4 is devoted to the proof of Theorem \ref{nocon2}. In Section 5, we will mention a few concluding remarks.
\section{Preliminaries}
We begin with the following fact.
\begin{fact}
Let $G$ be a  graph and $C$ be a longest cycle of length $c$. If $c<|V(G)|$, then $d_C(u)\le\lfloor c/2\rfloor$ for each $u\in V(G)\setminus V(C)$.
\end{fact}
 If $C$ is a longest cycle of length $c$ in $G$, then  let
 \[
 A_{\lfloor c/2\rfloor}=\{u\in V(G)\setminus V(C):d_C(u)=\lfloor c/2\rfloor\}.
 \]
Let us recall the following classical result by Bondy.
\begin{lemma}[Bondy \cite{Bondy}]\label{Bondy}
	Let $G$ be a graph on $n$ vertices and $C$ be a longest cycle of $G$ with order $c$. Then
	\[
	e(G-C)+e(G-C,C)\le\lfloor c/2 \rfloor(n-c).
	\]
\end{lemma}
Ma and Ning \cite{MN} proved a stability version of Lemma \ref{Bondy}. The first author, Ning, and the third author discovered the following variant for  2-connected graphs, see Lemma 7 in \cite{DNP}.
\begin{lemma}[Dou, Ning, Peng \cite{DNP}]\label{outcycle}
	Let $G$ be a $2$-connected graph on $n$ vertices and $C$ be a longest cycle in $G$ of length $c$ with $4 \leq c\le n-1$.
	
	\noindent
	$(1)$ If $A_{\lfloor c/2\rfloor}\neq\emptyset$, then
	\[e(G-C)+e(G-C,C)\le\lfloor c/2 \rfloor(n-c).\]
	\noindent
	$(2)$ If $A_{\lfloor c/2\rfloor}=\emptyset$, then
	\[e(G-C)+e(G-C,C)\le(\lfloor c/2 \rfloor-\tfrac{1}{2})(n-c).\]
\end{lemma}
We strengthen  Lemma \ref{outcycle} as follows.
\begin{lemma}\label{outcycle'}
	Let $G$ be a $2$-connected graph on $n$ vertices. If $C$ is a longest cycle in $G$ of length $c$ with $4 \leq c\le n-1$, then
	\[
	e(G-C)+e(G-C,C)\le \lfloor c/2 \rfloor|A_{\lfloor c/2\rfloor}|+(\lfloor c/2 \rfloor-\tfrac{1}{2})\left(n-c-|A_{\lfloor c/2\rfloor}|\right).
	\]
\end{lemma}
{\bf Proof:} If $A_{\lfloor c/2\rfloor}=\emptyset$, then we are done by Lemma \ref{outcycle}. For the case that  $A_{\lfloor c/2\rfloor}\neq\emptyset$, let $u\in A_{\lfloor c/2\rfloor}$ and $G_1=G-C$.   We  claim $N_{G_1}(u)=0$. Otherwise,  let $v\in V(G_1)$ be a neighbor of $u$.  As $G$ is $2$-connected,  there are two internally disjoint paths $P_1$ and $P_2$ from $v$ to $C$. Thus there exists a path, say $P_1$, avoiding $u$. We may assume that $\{y\}= P_1\cap C$. The assumption $u \in A_{\lfloor c/2\rfloor}$ implies that there is $x\in N_C(u)$ such that  $d_C(x,y) \in  \{1,2\}$, where $d_C(x,y)$ is the distance  between $x$ and $y$ on $C$. For $d_C(x,y)=1$, note that $$(C-xy)\cup(xuvP_1y)$$  is a longer cycle with length at least $c+1$.  For $d_C(x,y)= 2$, let $xwy$ be the path of length $2$ in $C$. Similarly,
$$(C-xwy)\cup(xuvP_1y)$$
  is a longer cycle of length at least $c+1$, a contradiction. Thus $N_{G_1}(u)=0$.

 Let $G'=G-A_{\lfloor c/2\rfloor}$. We know $G'$ is also 2-connected by $N_{G_1}(u)=0$. It is obvious that $C$ is a longest cycle in $G'$. If $G'$ is Hamiltonian, then $V(G-C)=A_{\lfloor c/2\rfloor}$. The conclusion is true. Thus we can assume $G'$ is not Hamiltonian. If we define
 $$A'_{\lfloor c/2\rfloor}=\{u\in V(G')\setminus V(C): d_C(u)=\lfloor c/2\rfloor\},$$ then  $A'_{\lfloor c/2\rfloor}=\emptyset$.   Lemma \ref{outcycle} yields that
 \begin{align*}
 	e(G'-C)+e(G'-C,C)&\le (\lfloor c/2 \rfloor-\tfrac{1}{2})\left(|V(G')|-c\right)\\
 	&=  (\lfloor c/2 \rfloor-\tfrac{1}{2})\left(n-c-|A_{\lfloor c/2\rfloor}|\right).
 \end{align*}
Therefore,
\begin{align*}
e(G-C)+e(G-C,C)&=\lfloor c/2 \rfloor|A_{\lfloor c/2\rfloor}|+e(G'-C)+e(G'-C,C)\\
&\le \lfloor c/2 \rfloor|A_{\lfloor c/2\rfloor}| + (\lfloor c/2 \rfloor-\tfrac{1}{2})\left(n-c-|A_{\lfloor c/2\rfloor}|\right).
\end{align*}
The lemma is proved. \hfill$\square$

Given a graph $F$, recall that
$$\h=\{F[S]: S \text{ is a  vertex covering of } F\}$$
 and
 $$\h'=\{F-\{u,v\}: u \textrm{ is adjacent to } v \textrm{ in } F\}.$$
If $G$ is an $F$-free graph and $C$ is a longest cycle of length $c$ in $G$,  then we shall estimate   the number of edges in $C$ as follows.
\begin{lemma}\label{cycle edge}
Let $G$ be an $F$-free graph on  $n$ vertices and $C$ be a longest cycle  of length $3 \leq c \leq n-1$ in $G$, where $p(F) \geq \lfloor c/2 \rfloor +1$.
	
	\noindent
	$(1)$  If $c$ is even and $|A_{c/2 }|\geq |V(F)|$, then
	$$e(G[C])\le \ex(c/2 ,\h)+(c/2)^2.$$
	\noindent	
	$(2)$ If $c$ is odd and $|A_{\lfloor c/2\rfloor}| \geq c|V(F)|$, then
	$$e(G[C])\le\ex(\lfloor c/2 \rfloor,\h)+\lfloor c/2 \rfloor\lceil c/2 \rceil+1.$$
	Moreover, if every graph in $\EX(\lfloor c/2 \rfloor,\h)$ contains at least one graph in $\h'$ as a subgraph, then the inequality is strict.
\end{lemma}
{\bf Proof:} (1) If $c$ is even, then  let $C=a_1,b_1,a_2,b_2,\ldots, a_{c/2},b_{c/2},a_1$.  Assume that $u$ is a vertex from $A_{c/2}$. Since $u$ can not be adjacent to consecutive vertices on $C$,   we  assume $L= N_C(u)=\{a_1,\ldots,a_i,\ldots,a_{c/2}\}$. Let $R=V(C)\setminus L$. Then $e(R)=0$. Otherwise, $G$ contains a longer cycle by including $u$, a contradiction. The assumption $|A_{ c/2}| \geq |V(F)|$ implies that $G[L]$ is  $\h$-free. Otherwise, $G[L]$ contains $F[S]$ as a subgraph for some vertex covering $S$.  Note that $F \setminus S$ is an independent set in $F$. As $A_{c/2}$ and  $L$ forms a complete bipartite graph, $G[A_{c/2} \cup L]$ will contain $F$ as a subgraph which is a contradiction to the assumption for $G$.
 It follows that
\begin{align*}
	e(G[C])&=e(L)+e(R)+e(L,R)\\
	&\le\ex(c/2,\h)+0+(c/2)^2\\
	&= \ex(c/2,\h)+(c/2)^2.
\end{align*}
$(2)$ We next consider the case that $c$ is odd.  It is trivial for $c=3$.  If $c\ge 5$, then  we  define the longest cycle $C$ as $v_1,v_2,\ldots,v_{c}$ with vertices ordered clockwisely.  For each $u \in A_{\lfloor c/2\rfloor}$, since $u$ can not be adjacent to consecutive vertices on $C$,   we may assume $N_C(u)=\{v_i,v_{i+3},v_{i+5},\ldots,v_{i+c-2}\}$ for some $1 \leq i \leq c$, where the addition is under modulo $c$. For each $1 \leq i \leq c$, let $A_i=\{u \in A_{\lfloor c/2\rfloor}: N_C(u)=\{v_i,v_{i+3},v_{i+5},\ldots,v_{i+c-2}\} \}$. Apparently, $A_{\lfloor c/2\rfloor}=\cup_{i=1}^c A_i$. Without loss of generality,  we may assume $|A_1| \geq |V(F)|$ because of the assumption $|A_{\lfloor c/2\rfloor}| \geq c|V(F)|$.
Similar to the previous case, let $L=\{v_1,v_{4},v_{6},\ldots,v_{c-1}\}$ and $R=V(C)\setminus L.$  Thus $G[L]$ is $\h$-free and $G[R]$  contains only an edge $v_2v_3$. We have
 \begin{align*}
 	e(G[C])&=e(L)+e(L,R)+e(R)\\
 	&\le\ex(\lfloor c/2\rfloor,\h)+\lfloor c/2 \rfloor\lceil c/2 \rceil+1.
 \end{align*}
 The equality holds only if $v_2v_3$ is an edge, $G[L]\in \EX(\lfloor c/2 \rfloor,\h)$ as well as $L$ and $R$ form a complete bipartite graph. If every graph in $\EX(\lfloor c/2 \rfloor,\h)$ contains  a member in $\h'$ as a subgraph, then the inequality is strict. Otherwise, the subgraph of $G$ induced by $L \cup \{v_2,v_3\}$ contains $F$ as a subgraph,  a contradiction.
 The lemma is proved. \hfill $\square$
%
%

To simplify  the proofs of Theorem \ref{odd} and Theorem \ref{even}, we need the following Lemma. 
\begin{lemma}\label{edgeG}
 Let $G$ be a  $2$-connected $n$-vertex graph which is $\C_{\geq k}$-free. Let $C$  be a longest cycle of length $c$ in $G$ and $|A_{\lfloor c/2\rfloor}|=m$.
	\item [(1)] If $c=k-1$, then \[e(G)\le e(G[C])+t\cdot m+(t-1/2)\left(n-c-m\right).\]
	\item [(2)] If $c=k-2$ and $k$ is even, then
	\[e(G)\le e(G[C])+t\cdot m+(t-1/2)\left(n-c-m\right).\]
\item [(3)] If  either $c\le k-3$ or  $c=k-2$ with $k$ being odd,  then
		\[e(G)\le (t-1/2)(n-1).\]
\end{lemma}
\noindent
{\bf Proof:} Part (1) follows from Lemma \ref{outcycle'} directly. For part (2), notice that $t=(k-2)/2$ and it also follows from  Lemma \ref{outcycle'}. Part (3) can be verified by  Theorem \ref{EGcycle}.   \hfill $\square$
\section{Proofs of Theorem \ref{odd} and  Theorem \ref{even}}
Fix $t=\lfloor (k-1)/2 \rfloor$. Let $\f=\{\C_{\ge k}, F\}$ with $p(F)\ge t+1$.  We first introduce the following construction of graphs. For each $T \in \EX(t,\h)$, let $G=T\vee I_{n-t}$.
Observe that $G$ is an $\f$-free 2-connected graph and $e(G)=\ex(t,\h)+t(n-t)$. Therefore, we need only to establish  upper bounds in Theorem \ref{odd} and Theorem \ref{even}.

  \noindent
{\bf Proof of Theorem \ref{odd}:} Notice that $k$ is odd. Let $G$ be a 2-connected $n$-vertex  graph which is $\f$-free.
A longest cycle of length $c$ in $G$ is denoted by $C$.  Assume that $|A_{\lfloor c/2\rfloor}|=m$.

We first consider the case that  $c=k-1$.
 If $m< c|V(F)|$, then by Lemma \ref{edgeG},  we have
\begin{align*}
e(G)&\le e(G[C])+t\cdot m+(t-1/2)\left(n-k+1-m\right)\\
    &\leq e(G[C])+ (t-1/2)n+m/2 \\
    &<\ex(t,\h)+t(n-t)
\end{align*}
for sufficiently large $n$. If $m \geq c|V(F)|$, then we show the following inequality by  combining  Lemma \ref{cycle edge} and Lemma \ref{edgeG}. Note that  $c=k-1$ is even and  $t=(k-1)/2$. Thus
\begin{equation}\label{ieq0}
\begin{aligned}
	e(G)&\le e(G[C])+t\cdot m+(t-1/2)\left(n-k+1-m\right)\\
	&\le\ex(t,\h)+t^2 +t\cdot m+(t-1/2)\left(n-k+1-m\right)\\
	&\leq \ex(t,\h)+t^2 +t\cdot m+(t-1/2)\left(n-2t-m\right)\\
	&\le\ex(t,\h)+t^2 +t\cdot m+t(n-2t-m)\\
	&=\ex(t,\h)+t(n-t).
\end{aligned}
\end{equation}

In the remaining case, by Lemma \ref{edgeG}, we have
\[e(G)\le (t-1/2)(n-1)\le\ex(t,\h)+t(n-t)\]
for sufficiently large $n$.
Thus
\[
\ex_{2\textrm{-conn}}(n,\f) = \ex(t,\h)+t(n-t),
\]
and the proof is complete. \hfill $\square$
\vspace{6pt}

\noindent
{\bf Proof of  Theorem \ref{even}:}
Notice that $k$ is even. Let $G$ be a 2-connected $n$-vertex graph which is $\f$-free. We assume that $C$ is a longest cycle of length $c$ in $G$ and $|A_{\lfloor c/2\rfloor}|=m$.
For the case where $c=k-1$, we can repeat the argument  in the proof of Theorem \ref{odd} to show the desired upper bound, here $t=(k-2)/2$ and  $e(C)\le\ex(t,\h)+t(k-1-t)+1$ in inequality \eqref{ieq0}.


 For $c=k-2$, we notice that $t=(k-2)/2$. If $m<|V(F)|$, then by Lemma \ref{edgeG} we have
 \begin{align*}
 	e(G)&\le e(G[C])+t\cdot m+(t-1/2)\left(n-k+2-m\right)\\
 	& \leq (t-1/2)n+m/2 \\
    &<\ex(t,\h)+t(n-t)
 \end{align*}
 for sufficiently large $n$. For $m \geq |V(F)|$, we use Lemma \ref{cycle edge} and Lemma \ref{edgeG} to show the following inequality. Note that  $c=k-2$ is even and $t=(k-2)/2$. Thus
 \begin{align*}
 	e(G)&\le e(G[C])+t\cdot m+(t-1/2)\left(n-k+2-m\right)\\
 	&\le\ex(t,\h)+t^2 +t\cdot m+(t-1/2)\left(n-k+2-m\right)\\
 	&\leq \ex(t,\h)+t^2 +t\cdot m+(t-1/2)\left(n-2t-m\right)\\
 	&\leq \leq \ex(t,\h)+t^2 +t\cdot m+t\left(n-2t-m\right)\\
 	&=\ex(t,\h)+t(n-t).
 \end{align*}

 For the leftover case, by Lemma \ref{edgeG}, we have
 \[e(G)\le (t-1/2)(n-1)\le\ex(t,\h)+t(n-t)\]
 for sufficiently large $n$.
  Thus
  \[
  \ex_{2\textrm{-conn}}(n,\f)\leq \ex(t,\h)+t(n-t)
  \]
for $c \leq k-2$. For $c=k-1$, we already showed
 \[
  \ex_{2\textrm{-conn}}(n,\f)\leq \ex(t,\h)+t(n-t)+1.
  \]
 If every graph in $\EX(t,\h)$ contains at least a member in $\h'$,   then the equality is strict in Lemma \ref{cycle edge} for $c=k-1$ being odd. It follows that
  \[e(G[C])\le\ex(t,\h)+t(k-1-t)\]
and
  \[e(G)\le\ex(t,\h)+t(n-t)\]
  for sufficiently large $n$ and thus $\ell=0$ in this case. The proof is complete.
\hfill $\square$
 \vspace{6pt}

%
%
 \section{Proof of Theorem \ref{nocon2}}
 Throughout the proof, we assume $\f=\{\C_{\ge k}, F\}$, where $F$ is 2-connected and $p(F)\ge t+1$. We again fix $t=\lfloor(k-1)/2\rfloor=\tfrac{t-2}{2}$.

 For the case where $k=4$,  Theorem \ref{EGcycle} gives us $\ex(n,\f) \leq \tfrac{3}{2}(n-1)$. Note that a friendship graph is $\f$-free. Thus $\ex(n,\f) \geq  \tfrac{3}{2}(n-1)$ for odd $n$ and $\ex(n,\f) \geq  \tfrac{3}{2}(n-2)$ for even $n$. The theorem follows in this case. For the rest of the proof, we assume $k \geq 6$.

 For the lower bound, it is easy to see that $\ex(n,\f)\ge \ex_{2\textrm{-conn}}(n,\f)=t\cdot n+O_k(1)$  by  Theorem \ref{even} for sufficient large $n$. Moreover, let $n-1=p(k-2)+q$, where $0\le q \le k-3$. Assume that $F_1\in\EX(k-1,F)$ and  $F_2\in\EX(q+1,F)$,   we know the graph that consists of $p$ copies of $F_1$ and a copy of $F_2$ sharing a common vertex is $\f$-free. This is true because $F$ is 2-connected. It follows that
 \begin{align*}
 	\ex(n,\f)&\ge p\cdot \ex(k-1,F)+\ex(q+1,F)\\
 	&=\frac{\ex(k-1,F)}{k-2}n+O_k(1).
 \end{align*}
 Thus $\ex(n,\f)\ge n\max\left\{\tfrac{k-2}{2},\frac{\ex(k-1,F)}{k-2}\right\}+O_k(1).$

  For the upper bound, let $G$ be an $n$-vertex graph which is $\f$-free with $\ex(n,\f)$ edges. As both $F$ and $\C_{\geq k}$ are 2-connected,  we may assume  $G$ is connected.  Let $B_1,\ldots,B_w$ be the blocks of $G$ and $|B_i|=b_i$. We next show
  \[
  e(G)\le(n-1)\cdot\max\left\{\frac{k-2}{2},\frac{\ex(k-1,H)}{k-2}\right\}.
  \]
%

%
For even  $k \geq 6$, recall that $\lfloor \tfrac{k-1}{2} \rfloor=\tfrac{k-2}{2}$.
If $b_i<k-1$, then $e(B_i)\le\frac{(b_i-1)b_i}{2}\le\frac{k-2}{2}(b_i-1)$. If $b_i=k-1$, then $e(B_i)\le\ex(k-1,F)$. For the remaining case that $b_i\ge k$, let $C$  be a longest cycle  of length $c$ in  $B_i$, where $c\le k-1$.

If $c\le k-2$, the Erd\H{o}s and Gallai Theorem implies that $e(B_i) \leq \tfrac{k-2}{2}(b_i-1)$. For  $c=k-1$, we combine the fact $e(B_i[C])\le \ex(k-1,F)$ and Lemma \ref{outcycle'} to get
\begin{align*}
	e(B_i)&\le e(B_i[C])+\lfloor c/2 \rfloor(b_i-k+1)\\
	&\le \ex(k-1,F)+\frac{k-2}{2}(b_i-k+1)\\
	&=\ex(k-1,F)+\frac{k-2}{2}(b_i-1)-\frac{(k-2)^2}{2}.
\end{align*}
Let \[A_1=\{B_i : b_i<k-1\},A_2=\{B_i : b_i=k-1\},\]
\[A_3=\{B_i : b_i\ge k \textrm{ and } c\le k-2\},A_4=\{ B_i: b_i\ge k \textrm{ and } c= k-1\}.\]
Apparently, $\sum_{i=1}^4|A_i|=w$. For $1 \leq i \leq 4$, let
 \[n_i=\sum_{B_j \in A_i }b_j \textrm{ and } q_i=|A_i|.\]
By the induction on the number of blocks, we can show
 \[
 \sum_{i=1}^{4}(n_i-q_i)=\sum_{i=1}^w (b_i-1)=n-1.
 \]
 Moreover, we claim that $n_4-q_4> q_4 (k-2).$ Actually,
 $$n_4-q_4=\sum_{i \in A_4}(b_i-1)> q_4(k-2)$$
  as $b_i\ge k$ for each $B_i\in A_4$.  By the same argument, we get that $q_2=\frac{n_2-q_2}{k-2}$. Now, notice that $$e(G)=\sum_{i=1}^4 \sum_{B_j \in A_i} e(B_j).$$
We  estimate each summand one by one as follows. For $i=1$, we have
\[
\sum_{B_j \in A_1}e(B_j) \leq \sum_{B_j \in A_1} \frac{k-2}{2}(b_j-1) =\frac{k-2}{2}(n_1-q_1).
\]
Similarly,
\[
\sum_{B_j \in A_2} e(B_j) \leq q_2\cdot\ex(k-1,F) \textrm{ and } \sum_{B_j \in A_3}e(B_j) \leq \frac{k-2}{2}(n_3-q_3).
\]
If $i=4$, then
\[
\sum_{B_j \in A_4} e(B_j) \leq \frac{k-2}{2}(n_4-q_4)+q_4\left(\ex(k-1,F)-\frac{(k-2)^2}{2}\right).
\]
Putting together, we obtain
\[
 e(G)
 \le \frac{k-2}{2}\left(\sum_{i=1}^{4}(n_i-q_i)-(n_2-q_2)\right)+\frac{n_2-q_2}{k-2}\ex(k-1,F)+q_4\left(\ex(k-1,F)-\frac{(k-2)^2}{2}\right).
\]
For the  inequality above, if   $\ex(k-1,F)-(k-2)^2/2\le 0$, then we have
\begin{align*}
		e(G)&\le\frac{k-2}{2}\left(\sum_{i=1}^{4}(n_i-q_i)-(n_2-q_2)\right)+\frac{n_2-q_2}{k-2}\ex(k-1,F) \\
            & \leq \left(\sum_{i=1}^{4}(n_i-q_i)\right) \cdot \max\left\{\frac{k-2}{2},\frac{\ex(k-1,F)}{k-2}\right\}\\
		&=(n-1)\cdot\max\left\{\frac{k-2}{2},\frac{\ex(k-1,F)}{k-2}\right\}.		
\end{align*}
The last equality holds as we already showed  $\sum_{i=1}^{4}(n_i-q_i)=n-1$. If  $\ex(k-1,H)-(k-2)^2/2> 0$, then we notice that $q_4<\frac{n_4-q_4}{k-2}$  as $n_4-q_4> q_4(k-2)$. Therefore,
\begin{align*}\label{ieq3} e(G)&\le\frac{k-2}{2}\left(\sum_{i=1}^{4}(n_i-q_i)-(n_2-q_2)\right)+\frac{n_2-q_2}{k-2}\ex(k-1,F)+\frac{n_4-q_4}{k-2}\left(\ex(k-1,F)-\frac{(k-2)^2}{2}\right)\nonumber\\
	&= \frac{k-2}{2}(n_1-q_1+n_3-q_3)+ \frac{\ex(k-1,F)}{k-2}(n_2-q_2+n_4-q_4)\\
    &\leq \left(\sum_{i=1}^{4}(n_i-q_i)\right) \cdot \max\left\{\frac{k-2}{2},\frac{\ex(k-1,F)}{k-2}\right\}\\
	&=(n-1)\cdot\max\left\{\frac{k-2}{2},\frac{\ex(k-1,F)}{k-2}\right\}.
\end{align*}
The upper bound is proved. \hfill$\square$
\section{Concluding remarks}
In this paper, we proved results on the Tur\'an number of $\{\C_{\geq k},F\}$ for a general 2-connected graph $F$.
For Theorem \ref{nocon2}, the assumption $F$ being 2-connected  is only needed in the proof of the lower bound. If one only concerns the upper bound for $\ex(n, \{\C_{\geq k},F\})$, then we can drop the assumption $F$ being 2-connected and  show
\[
 \ex(n, \{\C_{\geq k},F\}) \leq n\max\left\{\frac{k-2}{2},\frac{\ex(k-1,F)}{k-2}\right\}
\]
for sufficiently large $n$. If $F$ is not 2-connected, then Theorem \ref{EGcycle} gives a lower bound
\[
\ex(n, \{\C_{\geq k},F\}) \geq \frac{k-2}{2}n+O_k(1).
\]
If further $\tfrac{\ex(k-1,F)}{k-2} \leq \tfrac{k-2}{2}$, then
\[
\ex(n, \{\C_{\geq k},F\}) = \frac{k-2}{2}n+O_k(1).
\]
Therefore, it is an interesting question to determine $\ex(n, \{\C_{\geq k},F\})$ up to an additive constant term, where $F$ is not 2-connected, $\tfrac{\ex(k-1,F)}{k-2} > \tfrac{k-2}{2}$, and $p(F) \geq \lfloor (k-1)/2 \rfloor+1$.

In the study of Tur\'an number of paths and cycles, it is  well known  that one can deduce $\ex(n,P_k)$ from $\ex(n, \C_{\geq k+1})$.  As Liu and Kang \cite{LK} already showed
\[
\ex(n,\{P_k,F\})=n\max\left\{\left \lfloor \frac{k}{2} \right \rfloor-1,\frac{\ex(k-1,F)}{k-1}\right\}+O_k(1)
\]
for sufficiently large $n$, where $F$ is a connected graph with $p(F) \geq \lfloor \tfrac{k}{2}  \rfloor$. A natural question is  whether one can deduce $\ex(n,\{P_k,F\})$ from $\ex(n,\{\C_{\geq k+1},F\})$. In \cite{DNP}, the authors gave a positive answer to the question above for $F=K_r$ and $r \geq \lfloor k/2 \rfloor+1$. In general, for a graph $F$, let $F'$ be  $F \vee \{v\} $ which is the join of $F$ and a single vertex. If $F$ is connected and $p(F) \geq \lfloor \tfrac{k}{2}  \rfloor$, then $F'$ is 2-connected and $p(F') \geq \lfloor \tfrac{k}{2}  \rfloor+1$. If  $F$ further satisfies $$\ex(k,F')=\ex(k-1,F)+k-1,$$
 then $\ex(n,\{P_k,F\})$ is still a consequence of $\ex(n,\{\C_{\geq k+1},F'\})$ by using the folklore trick. Note that $K_r$ with $r \geq \lfloor k/2 \rfloor+1$ indeed satisfies the condition above. A challenging question in this direction is whether it is still the case for a general graph $F$.

\end{document}